\PassOptionsToPackage{hyphens}{url}
\documentclass[a4paper,USenglish,cleveref,autoref,thm-restate,nolineno]{socg-lipics-v2021}

\pdfoutput=1 
\hideLIPIcs  

\title{Bounded Twin-Width Refinement of Triangulations} 

\author{\'Edouard Bonnet}{Univ Lyon, CNRS, ENS de Lyon, Université Claude Bernard Lyon 1, LIP UMR5668, France \and \url{http://perso.ens-lyon.fr/edouard.bonnet/}}{edouard.bonnet@ens-lyon.fr}{https://orcid.org/0000-0002-1653-5822}{The author has been supported by the French National Research Agency through the project TWIN-WIDTH with reference number ANR-21-CE48-0014.}

\author{Kristóf Huszár}{Institute of Geometry, Graz University of Technology, Austria \and \url{https://kristofhuszar.github.io/}}{kristof.huszar@tugraz.at}{https://orcid.org/0000-0002-5445-5057}{}
\authorrunning{\'E. Bonnet and K. Huszár} 

\Copyright{\'Edouard Bonnet and Kristóf Huszár} 



\category{} 

\relatedversion{} 

\nolinenumbers 

\usepackage{mathtools}
\usepackage{overpic}
\usepackage{enumitem}

\newcommand{\bigsimplex}{\Sigma}

\newcommand{\dual}{\Gamma}

\newcommand{\manifold}{\mathscr{M}}

\newcommand{\orig}{\mathscr{O}}

\newcommand{\pillow}{\nabla}

\newcommand{\redpath}{\mathfrak{P}}
\newcommand{\refedge}{\mathfrak{R}}
\newcommand{\seq}{\mathfrak{S}}

\newcommand{\simplex}{\sigma}

\newcommand{\tri}{\mathscr{T}}
\newcommand{\tw}{\operatorname{tw}}
\newcommand{\tww}{\operatorname{tww}}

\newcommand{\width}{\operatorname{w}}

\newcommand{\mynum}[1]{{\color{lipicsGray}\sffamily\bfseries{#1}}}

\renewcommand{\leq}{\leqslant}
\renewcommand{\geq}{\geqslant}

\newcommand\smallscale{0.884615384615} 

\begin{document}

\maketitle

\begin{abstract}
We prove that any triangulation $\tri$ with size $|\tri|$ of a $d$-dimensional manifold can be refined to a triangulation $\tri^\ast$ of the same manifold, with size $|\tri^\ast| \leq 2d \cdot |\tri| \cdot \log_2 |\tri|$ and dual graph $\dual(\tri^\ast)$ of twin-width at most four. The construction is based on work of Bergé, Bonnet, and Déprés, who showed that any $(\geq 2 \log_2 n)$-subdivision of an $n$-vertex graph has twin-width at most four.
\end{abstract}

\section{Introduction}
\label{sec:intro}

Triangulations play a fundamental role in the study of manifolds. This is particularly true in three dimensions, where every topological 3-manifold admits a triangulation \cite{moise1952affine},\footnote{Topological manifolds of dimension at least four cannot always be triangulated \cite{manolescu2016lectures}.} and where a substantial body of work over the past decades has been guided by computer experiments (see \cite{basilio2024totally,costantino2025learning,futer2025excluding,he2025practical,kegel2026volume,li2026complete} for recent examples; also on 4-manifolds \cite{burke2025practical,burke2026small}), aided by specialized software \cite{Regina,SnapPy}, cf.\ \cite{computop}. In this context, triangulations provide a useful discrete description of a topological space. However, there is a caveat: every triangulable manifold has infinitely many combinatorially distinct triangulations, and different triangulations of the same manifold can exhibit strikingly different algorithmic behavior. While some triangulations may allow us to quickly obtain topological information about the underlying manifold, others may be so complicated that they are prohibitively difficult to explore.\footnote{See \cite{MFOReports-Computational,triangulations-book,MFOReports-Triangulations} for a panoramic overview of how triangulations are used in computational settings.}

Therefore, it is often desirable to represent manifolds by triangulations that are as simple as possible, where simplicity is measured not primarily by size,\footnote{The \emph{size} of a triangulation is defined as the number of its top-dimensional simplices.} but by structural graph parameters. Indeed, several \emph{fixed-parameter tractable} (FPT) algorithms are now known that efficiently compute topological invariants from triangulations of bounded \emph{treewidth},\footnote{Treewidth is a graph parameter that quantifies the similarity of a given graph to any tree. The \emph{treewidth of a triangulation} $\tri$ is defined as the treewidth of its \emph{dual graph} $\dual(\tri)$. See \Cref{sec:prelims} for details.} even if the corresponding problems are computationally hard in general \cite{bagchi2016efficient,burton2017courcelle, burton2016parameterized, burton2018algorithms, pettersson2014fixed, burton2013complexity,huszar2026sparse,philip2026eth}.

Finding triangulations with small treewidth, however, is challenging. While several results assert the existence of triangulations whose treewidth is bounded above by certain topological invariants of the underlying $3$-manifold, such as the Heegaard genus \cite{he2024algorithm,huszar2019manifold} or the hyperbolic volume \cite{huszar2022pathwidth,maria2019treewidth}, there are complementary results showing that, under certain conditions, a $3$-manifold with complicated topology cannot admit a triangulation with small treewidth \cite{huszar2025width,huszar2019treewidth}. More importantly, in a computational setting, a manifold is typically given by a triangulation (or by other discrete data structures, cf.\ \cite{cannon2003heegaard,ennes2026compressed}), thus the actual task would be to \emph{retriangulate} the manifold so as to reduce its treewidth before applying any of the aforementioned FPT algorithms. However, to the best of our knowledge, there is currently no exact or heuristic method that can either reduce the treewidth of a given triangulation or certify that no such retriangulation is possible.

In this work, we focus on another structural graph parameter, namely \emph{twin-width}. Introduced in \cite{bonnet2022twin-width-1}, twin-width has attracted considerable attention and has since been studied in a wide range of settings (see \cite{bonnet2024habil} for a survey). Motivated by the considerations discussed above, it was recently shown that every $d$-dimensional \emph{smooth} manifold has a triangulation with twin-width $d^{O(d)}$ \cite[Theorem~1]{bonnet2025smooth}.
Thus, every $3$-manifold has a triangulation with twin-width bounded above by a universal constant \cite[Corollary~2]{bonnet2025smooth}, in sharp contrast to the case of treewidth, cf.\ \cite{huszar2025width,huszar2019treewidth}.
This bound in \cite{bonnet2025smooth} was derived from an analysis of Whitney's classical method for triangulating smooth manifolds \cite{whitney1957geometric} (cf.\ \cite{boissonnat2021triangulating}). Now, Whitney's method assumes that the $d$-manifold to be triangulated is already smoothly embedded in $\mathbb{R}^{2d}$. Although such an embedding always exists for smooth manifolds, in a setting where the manifold is given by a triangulation, finding a smooth embedding and subsequently carrying out Whitney's construction may be impractical, or even computationally infeasible.

Here we show that any triangulation of a $d$-manifold can directly be refined into another triangulation of the same manifold whose twin-width is at most four. First, let us recall that the \emph{size} $|\tri|$ of a $d$-manifold triangulation $\tri$ is defined as the number of its $d$-simplices, while the twin-width of $\tri$ is defined as the twin-width of its dual graph $\dual(\tri)$, denoted $\tww(\dual(\tri))$, cf.\ \Cref{sec:prelims} for precise definitions. With this, our main result is the following theorem.

\begin{theorem}
\label{thm:const-tww}
Any triangulation $\tri$ of a $d$-dimensional manifold $\manifold$ with size $|\tri|$ can be refined to a triangulation $\tri^\ast$ of $\manifold$ of size $|\tri^\ast| = O(d \cdot |\tri| \cdot \log_2 |\tri|)$, such that $\tww(\dual(\tri^\ast)) \leq 4$.
\end{theorem}

\subparagraph*{Outline of the proof} The proof of \Cref{thm:const-tww} is based on a simple construction inspired by \cite[Theorem 2]{berge2022deciding}, which states that a $(\geq 2 \log_2 n)$-subdivision of any $n$-vertex graph has twin-width at most four. Given a triangulation $\tri$ of a $d$-manifold $\manifold$, we construct another triangulation $\tri^\ast$ of $\manifold$ as follows. For each face-gluing $\varphi$ of $\tri$, we first ``unglue'' the interiors of the two $(d-1)$-simplices identified by $\varphi$, and then fill the resulting void by appropriately gluing in a complex $\pillow_{s,d}$, which we call an \emph{$s$-stacked $d$-pillow}, where $s$ is sufficiently large. This yields a triangulation $\tri^\ast$, whose dual graph $\dual(\tri^\ast)$ is reminiscent of $\dual(\tri)$ with its edges subdivided $s$ times. We show $\tww(\dual(\tri^\ast)) \leq 4$ by exhibiting a contraction sequence analogous to the one used in the proof of \cite[Theorem 2]{berge2022deciding}. Bounding the size of $\tri^\ast$ is immediate.

\subparagraph*{Discussion of the result} The upshot of \Cref{thm:const-tww} is that it provides a direct construction of a constant-twin-width triangulation $\tri^\ast$ of any triangulable $d$-manifold $\manifold$, given an arbitrary triangulation $\tri$ of $\manifold$. Moreover, for every fixed $d$, the size of $\tri^\ast$ is polynomially bounded in that of $\tri$. In particular, the construction can be carried out effectively in polynomial time.

At the same time, \Cref{thm:const-tww} may also be viewed as a potentially negative result from the perspective of parameterized algorithms: such a construction provides strong evidence against the possibility of FPT algorithms parameterized by twin-width for computing ``powerful'' topological invariants of manifolds given by triangulations, many of which are NP-hard to compute in general. Indeed, several such topological invariants can be expressed in monadic second-order logic, for which FPT algorithms parameterized by \emph{treewidth} are available \cite{burton2017courcelle}. Bounded twin-width, by contrast, is particularly well suited to first-order model checking \cite{bonnet2022twin-width-1}, and thus does not seem to offer the same algorithmic advantages for these invariants.

Thus, the main significance of \Cref{thm:const-tww} appears to lie in what it reveals about twin-width, rather than in immediate algorithmic consequences for computational topology.

\section{Preliminaries}
\label{sec:prelims}

Partially based on \cite[Sections 3 and 4]{bonnet2022twin-width-1}, we first review the notions we build upon.

\subsection{Graphs and triagraphs}
\label{ssec:graphs}

\subparagraph*{Graphs} A \emph{graph}, or, more specifically, a \emph{multigraph} $G$, is a pair $G = (V,E)$ consisting of a finite set $V$ of \emph{vertices} and a finite multiset $E$ of two-element submultisets of $V$, called \emph{edges}. An edge of the form $\{v,v\}$ is called a \emph{loop}, while an element of $E$ with multiplicity greater than one is called a \emph{multiedge}. A \emph{simple graph} is a graph with neither loops nor multiedges. The \emph{degree} $\deg(v)$ of a vertex $v \in V$ is the number of edges containing it, counted with multiplicity. Graphs in this paper are typically encountered as dual graphs of triangulations, see \Cref{ssec:tri} for details.

\subparagraph*{Trigraphs} Informally, a \emph{trigraph} $G$ is a (multi)graph together with a bicoloring of its edges. More precisely, it is a triple $G = (V,E,R)$, where, as before, $V$ is a finite set of \emph{vertices}, and $E$ and $R$ are two disjoint, finite multisets of two-element submultisets of $V$, called \emph{black edges} and \emph{red edges}, respectively.
The \emph{red degree} $\deg_R(v)$ of a vertex $v \in V$ equals the number of red edges incident to $v$, counted with multiplicity. Crucial to the definition of twin-width (see \Cref{ssec:tww}), we let $\deg_R(G) = \max_{v \in V(G)}{\deg_R(v)}$ denote the \emph{maximum red degree} of $G$.

\subsection{Contraction sequences and twin-width}
\label{ssec:tww}

Now, let $G = (V,E,R)$ be a trigraph and $u,v \in V$ be two arbitrary distinct vertices of $G$. We say that the trigraph $G / u,v = (V',E',R')$, called a \emph{contraction} of $G$, is obtained from $G$ by \emph{contracting} $u$ and $v$ into a new vertex $w$ if \mynum{1.}~$V' = (V \setminus \{u,v\}) \cup \{w\}$, \mynum{2.}~$G - \{u,v\} = (G/u,v) - \{w\}$ and \mynum{3.}~for any $x \in V'\setminus\{w\} = V\setminus\{u,v\}$ we have
\begin{itemize}
	\item $\{w,x\} \in E'$ if and only if $\{u,x\} \in E$ and $\{v,x\} \in E$,
	\item $\{w,x\} \notin E' \cup R'$ if and only if $\{u,x\} \notin E \cup R$ and $\{v,x\} \notin E \cup R$, and
	\item $\{w,x\} \in R'$ otherwise.
\end{itemize}

\begin{figure}[ht]
	\includegraphics[width=1\textwidth]{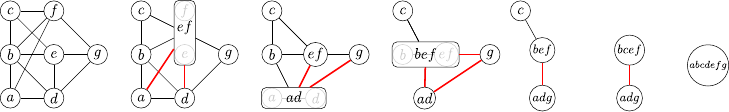}
	
	\medskip
	
	\caption{Left to right: a contraction sequence of  width two. Reproduced from \cite[p.~24]{bonnet2024habil}.}
	\label{fig:contraction}
\end{figure}

A sequence $\seq = (G_1,\ldots,G_m)$ of trigraphs is a \emph{contraction sequence} if $G_{i+1}$ is a contraction of $G_i$ for every $1 \leq i \leq m-1$. Note that $|V(G_{i+1})| = |V(G_i)|-1$.
The \emph{width} $\width(\seq)$ of a contraction sequence $\seq = (G_1,\ldots,G_m)$ is defined as $\width(\seq) = \max_{1 \leq i \leq m}\deg_R(G_i)$, i.e., the largest red degree of any vertex of any trigraph in $\seq$.
Now, the \emph{twin-width} $\tww(G)$ of a trigraph $G$ is the smallest width of any contraction sequence $(G_1,\ldots,G_{|V(G)|})$, where $G_1 = G$ and $G_{|V(G)|}$ consists of a single vertex. See \Cref{fig:contraction} for an example.

\subsection{Triangulations}
\label{ssec:tri}

\begin{figure}[htbp]
\begin{subfigure}{0.55\textwidth}
    \centering
    \includegraphics[scale=\smallscale]{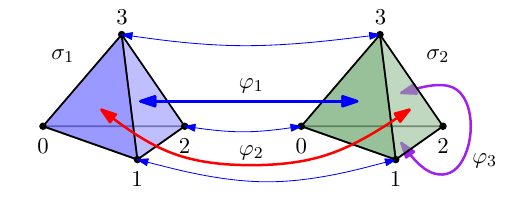}
    \caption{A triangulation $\tri$ with tetrahedra $\simplex_1$ and $\simplex_2$, and face gluing maps $\varphi_1$, $\varphi_2$ and $\varphi_3$ (exploded view).}
    \label{fig:tetrahedra}
\end{subfigure}%
\hfill%
\begin{subfigure}{0.4\textwidth}
    \centering
    \raisebox{0\height}{\includegraphics[scale=\smallscale]{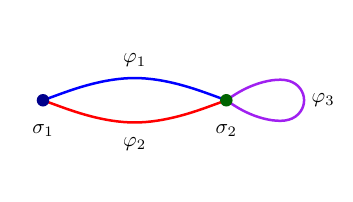}}
    \caption{The dual graph $\dual(\tri)$ of the triangulation $\tri$. $\dual(\tri)$ is generally a multigraph.}
    \label{fig:tetrahedra-dual}
\end{subfigure}%
\caption{A triangulation of a 3-manifold and its dual graph. Based on \cite[Figure 1]{huszar2022pathwidth}.}
\label{fig:gen-tri}
\end{figure}

Let $\bigsimplex = \{\simplex_1,\ldots,\simplex_n\}$ be a finite collection of vertex-labeled $d$-simplices and $\Phi = \{\varphi_1,\ldots,\varphi_m\}$ be a collection of $m \leq \frac{n(d+1)}{2}$ \emph{face gluing maps}, where each $\varphi_i$ is a simplicial isomorphism between two $(d-1)$-faces of the $d$-simplices in $\bigsimplex$. Each such $(d-1)$-face is the domain or the range of at most one $\varphi_i$, however, they do not need to be faces of distinct $d$-simplices. The quotient space $\tri = \bigsimplex/\Phi$ is called a \emph{triangulation}\footnote{Such triangulations are also called \emph{generalized}, \emph{semi-simplicial}, or \emph{singular} triangulations in the literature.} of the $d$-manifold $\manifold$ if its geometric realization $\|\tri\|$ is homeomorphic to $\manifold$, see \Cref{fig:tetrahedra} for an example.

Given a triangulation $\tri$, its \emph{dual graph} $\dual(\tri)$ is the multigraph with vertex set $\bigsimplex$ and edges corresponding to $\Phi$: two vertices are connected by an edge if and only if there is a face gluing map between the $(d-1)$-faces of the corresponding $d$-simplices, see \Cref{fig:tetrahedra-dual}.

\begin{remark*}
Triangulations are ubiquitous in computational topology. Note that a triangulation, as defined above, is not necessarily a simplicial complex. However, the second iterated barycentric subdivision $\tri''$ of any triangulation $\tri$ is always a simplicial complex. Triangulations in the present sense allow the description of manifolds using much fewer simplices than simplicial complexes, benefiting computations. For example, more than 13~thousand 3-manifolds can be described by a triangulation with 11 tetrahedra or less \cite[Theorem~9]{burton2011detecting}.
\end{remark*}

\section{\texorpdfstring{The proof of \Cref{thm:const-tww}}{The proof of the main result}}
\label{sec:const-tww}

We are ready to prove our main result, \Cref{thm:const-tww}, according to which any given triangulation of a $d$-dimensional manifold $\manifold$ can be refined into a \emph{triangulation} of $\manifold$, whose dual graph has twin-width at most four.
One of the main ingredients of the proof is the following theorem.
\begin{theorem}[{\cite[Theorem 2]{berge2022deciding}}]
\label{thm:log-subdiv}
Any $(\geq 2 \log_2 n)$-subdivision of an $n$-vertex graph has twin-width at most $4$.
\end{theorem}

\begin{proof}[Proof of \Cref{thm:const-tww}] 

The proof consists of two parts. First, we describe the refinement $\tri^\ast$ of $\tri$. Then we show that its dual graph $\dual(\tri^\ast)$ satisfies $\tww(\dual(\tri^\ast)) \leq 4.$

\subparagraph*{Part I. The construction of $\boldsymbol{\tri^*}$ from $\boldsymbol{\tri}$}

Let $\tri$ be a (simplicial or generalized) triangulation of a piecewise-linear $d$-manifold $\manifold$. We aim to refine $\tri$ into a triangulation $\tri^\ast$, so that we can use \Cref{thm:log-subdiv} to establish the upper bound $\tww(\dual(\tri^\ast)) \leq 4$. We achieve this refinement by ``inflating'' the triangulation $\tri$ at each of its face gluings with a \emph{$(2 \log_2 n)$-stacked $d$-pillow}, a triangulation of the $d$-ball we now describe.

\begin{description}[leftmargin=0pt]

	\item[$\boldsymbol{d}$-pillows] A \emph{$d$-pillow} is a triangulation $\pillow_d$ of the $d$-dimensional ball obtained by gluing together two $d$-simplices $\simplex_1$ and $\simplex_2$ as follows. Let $\{u_0,\ldots,u_d\}$ be the vertex set of $\simplex_1$ and $\{v_0,\ldots,v_d\}$ be that of $\simplex_2$. For $i \in \{1,2\}$ and $j \in \{1,\ldots,d\}$ let $\tau^j_i$ denote the $(d-1)$-face of $\simplex_i$ that does not contain the vertex $u_j$ (if $i=1$) or $v_j$ (if $i=2$). Let $\varphi_j \colon \tau_1^j \rightarrow \tau_2^j$ be the face gluing map for which $\varphi_j(u_k) = v_k$ for every $k \in \{0,\ldots,d\}\setminus \{j\}$. Letting $\bigsimplex = \{\simplex_1,\simplex_2\}$ and $\Phi = \{\varphi_1,\ldots,\varphi_d\}$, a \emph{$d$-pillow} is the triangulation $\pillow_d = \Sigma / \Phi$ (\Cref{fig:pillow}).
\begin{figure}[ht]
	\centering
	\begin{subfigure}{.45\textwidth}
		\centering
		\includegraphics[scale=\smallscale]{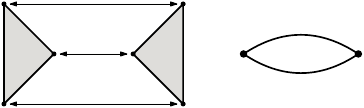}
		\smallskip
		\caption{A 2-pillow and its dual graph.}
	\end{subfigure}%
	\hfill%
	\begin{subfigure}{.5\textwidth}
		\centering
		\includegraphics[scale=\smallscale]{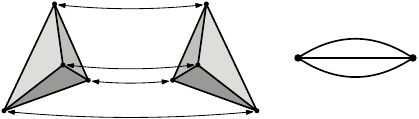}
		\smallskip
		\caption{A 3-pillow and its dual graph.}
	\end{subfigure}
	\caption{Exploded views of a 2-pillow and a 3-pillow along with their dual graphs.}
	\label{fig:pillow}
\end{figure}

\item[$\boldsymbol{s}$-stacked $\boldsymbol{d}$-pillows] By sequentially gluing them together along their boundary simplices in a natural way, $d$-pillows can be stacked on top of each other. More precisely, let $\pillow_d^{(1)},\ldots,\pillow_d^{(s)}$ be $d$-pillows, where $s \in \mathbb{N}$ is greater than 1. For each $j \in \{1,\ldots,s\}$, let $t_1^{(j)}$ and $t_2^{(j)}$ denote the two $(d-1)$-simplices that form the boundary of $\pillow_d^{(j)}$, and let $u_1^{(j)},\ldots,u_d^{(j)}$ denote the boundary vertices of $\pillow_d^{(j)}$. For $j \in \{1,\ldots,s-1\}$, let $\phi_j \colon t_2^{(j)} \rightarrow t_1^{(j+1)}$ be a face gluing map for which $\phi_j (u_k^{(j)}) = u_k^{(j+1)}$ for every $k \in \{1,\ldots,d\}$. Akin to the description of triangulations, letting $\pillow = \{\pillow_d^{(1)},\ldots,\pillow_d^{(s)}\}$ and $\Phi = \{\phi_j \colon 1 \leq j \leq s-1\}$, an \emph{$s$-stacked $d$-pillow} is the quotient space $\pillow_{s,d} = \pillow / \Phi$ (\Cref{fig:stacked-pillow}).
\begin{figure}[ht]
	\centering
	\includegraphics[scale=\smallscale]{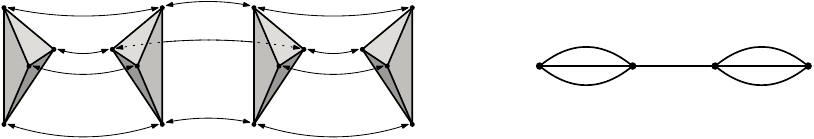}
	\caption{A 2-stacked 3-pillow (left; exploded view) and its dual graph (right).}
	\label{fig:stacked-pillow}
\end{figure}

\item[$\boldsymbol{s}$-refinements of triangulations] Given a triangulation $\tri = \bigsimplex / \Phi$ of a $d$-manifold $\manifold$, its \emph{$s$-refinement} $R_s(\tri)$ is another triangulation of $\manifold$ obtained as follows. For each face gluing map $\varphi \in \Phi$, which is a simplicial isomorphism $\varphi \colon \tau_1 \rightarrow \tau_2$ between two $(d-1)$-faces $\tau_1$ and $\tau_2$ of the $d$-simplices in $\tri$, take an $s$-stacked $d$-pillow $\pillow_{s,d}^{(\varphi)}$. Let $t_1$ and $t_2$ denote the two boundary $(d-1)$-simplices of $\pillow_{s,d}^{(\varphi)}$. Further, let $\tau_1(0) = \{u_1,\ldots,u_d\}$ and $t_1(0) = t_2(0) = \{v_1,\ldots,v_d\}$ be the vertices of the respective $(d-1)$-simplices. (As $t_1$ and $t_2$ have a shared boundary in $\pillow_{s,d}^{(\varphi)}$, we choose identical vertex labels for them.) Now consider two simplicial isomorphisms $\psi_\varphi^{(1)}\colon \tau_1 \rightarrow t_1$ and $\psi_\varphi^{(2)}\colon t_2 \rightarrow \tau_2$ defined as
\[
	\psi_\varphi^{(1)}(u_i) = v_i\quad \text{and} \quad \psi_\varphi^{(2)}(v_i) =  \varphi(u_i)\quad \text{for} \quad 1 \leq i \leq d.
\] 
Let $\pillow = \{\pillow_{s,d}^{(\varphi)} \colon \varphi \in \Phi\}$ and $\Psi = \bigcup_{\varphi \in \Phi}\{\psi_\varphi^{(1)},\psi_\varphi^{(2)}\}$. With this, the \emph{$s$-refinement} $R_s(\tri)$ of $\tri = \bigsimplex / \Phi$ is defined as the quotient space $R_s(\tri) = (\bigsimplex \cup \pillow)/\Psi$. See, e.g., \Cref{fig:refinement-2d}.

\begin{figure}[ht]
\centering
\begin{subfigure}{.32\textwidth}
	\centering
	\includegraphics[scale=\smallscale]{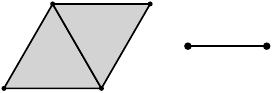}
	\smallskip
	\caption{A triangulation $\tri$ with two triangles and its dual graph $\dual(\tri)$.}
\end{subfigure}\hfill
\begin{subfigure}{.5\textwidth}
	\centering
	\includegraphics[scale=\smallscale]{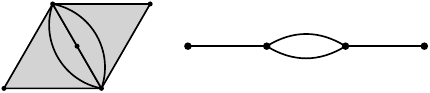}
	\smallskip
	\caption{The refinement $R_1(\tri)$ of $\tri$ via a 2-pillow at the original face gluing and its dual graph $\dual(R_1(\tri))$.}
\end{subfigure}

\bigskip

\begin{subfigure}{.8\textwidth}
	\centering
	\includegraphics[scale=\smallscale]{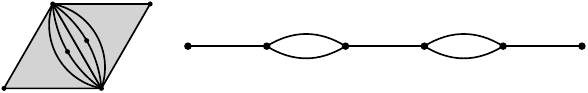}
	\smallskip
	\caption{The refinement $R_2(\tri)$ of $\tri$ via a 2-stacked 2-pillow and its dual graph $\dual(R_2(\tri))$.}
\end{subfigure}

\caption{A 2-dimensional triangulation $\tri$, its refinements $R_1(\tri)$, $R_2(\tri)$, and their dual graphs.}
\label{fig:refinement-2d}
\end{figure}

\end{description}

To conclude the construction, for a triangulation $\tri$ of a $d$-manifold $\manifold$ consisting of $n$ $d$-simplices we set $\tri^\ast$ to be the $(2\lceil \log_2 n \rceil-1)$-refinement\footnote{The choice of the value $2\lceil \log_2 n \rceil-1$ is motivated by the proof technique we rely on, see Part II.} of $\tri$. That is
\begin{align}
	\tri^\ast = R_{2\lceil \log_2 n \rceil - 1}(\tri).
\end{align}

\begin{description}[leftmargin=0pt]
	\item[Bounding the size of $\boldsymbol{\tri^\ast}$] Taking the $s$-refinement of a $d$-dimensional triangulation $\tri$ equals inserting an $s$-stacked $d$-pillow at each face gluing of $\tri$. Since a $d$-dimensional triangulation with $|\tri| = n$ has at most $(d+1)n/2$ face gluings and an $s$-stacked $d$-pillow consists of $2s$ $d$-simplices, we obtain that, for the number of $d$-simplices of $R_s(\tri)$ we have 
\[
	|R_s(\tri)| \leq (d+1)ns + n.
\]	
	In particular, for $\tri^\ast = R_{2\lceil \log_2 n \rceil - 1}(\tri)$ we have
\begin{align}
	|\tri^\ast| \leq (d+1)n(2\lceil \log_2 n \rceil - 1) + n = O(d n \log n).
\end{align}
\end{description}

\subparagraph*{Part II. Analysis of the twin-width of $\boldsymbol{\dual(\tri^*)}$} We show $\tw(\tri^\ast) \leq 4$ using the technique developed in \cite[Section~3]{berge2022deciding}, adapted to our setting. First we introduce some terminology.

\begin{description}[leftmargin=0pt]
	\item[Additional terminology] From the construction it is clear that, given a $d$-dimensional triangulation  $\tri$, the dual graph $\dual(R_s(\tri))$ of its $s$-re\-fine\-ment $R_s(\tri)$ is obtained from $\dual(\tri)$ by replacing each edge with a path of length $2s+1$, where every second edge is a multi-edge of multiplicity $d$, see, e.g., \Cref{fig:refined-edge}. We call such a path an \emph{$s$-refined edge}, or simply a \emph{refined edge}.
We call the vertices of $\dual(R_s(\tri))$ not incident to a multi-edge \emph{original vertices}, as they correspond to vertices of $\dual(\tri)$, see \Cref{fig:refined-edge,fig:refinement-example-2}.
Sequentially contracting\footnote{Here the term `contracting' is meant in the sense of \Cref{ssec:tww}.} each pair of vertices along the path that span a multi-edge yields a \emph{red path} of length $s+1$, i.e., a path with only red edges, cf.\ \Cref{fig:refined-edge-contracted}. We refer to this process as \emph{contacting the beads} along the given refined edge.
Although the dimension $d$ is not indicated in the notation, this shall not cause confusion. In particular, the proof is formally the same for every value of $d$.
\begin{figure}[ht]
\centering
\begin{subfigure}{.7\textwidth}
	\centering
	\includegraphics[scale=\smallscale]{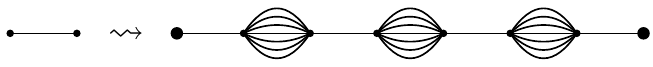}

	\smallskip
	\caption{Taking the $3$-refinement of a $7$-dimensional triangulation $\tri$ at an edge of its dual graph $\dual(\tri)$. Here the two original vertices are drawn larger.}
	\label{fig:refined-edge}
\end{subfigure}\hfill%
\begin{subfigure}{.25\textwidth}
	\centering
	\includegraphics[scale=\smallscale]{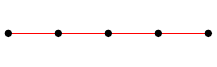}

	\smallskip
	\caption{A red path of length four (where $4 =s+1$).}
	\label{fig:refined-edge-contracted}
\end{subfigure}

\caption{$s$-refinement of an edge ($s=3$) and the red path obtained by contracting its beads.}
\end{figure}

\end{description}

\begin{figure}[ht]
\centering
\begin{subfigure}[t]{.25\textwidth}
	\centering
	\includegraphics[scale=\smallscale]{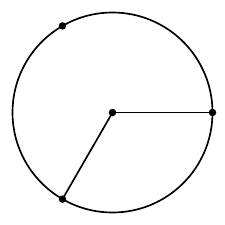}

	\smallskip
	\caption{The dual graph of a triangulated surface $\tri$.}
\end{subfigure}\hfill%
\begin{subfigure}[t]{.27\textwidth}
	\centering
	\includegraphics[scale=\smallscale]{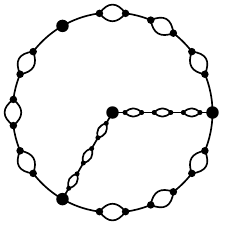}

	\smallskip
	\caption{The dual graph $\dual(R_3(\tri))$ of the $3$-refinement of $\tri$.}
	\label{fig:refinement-example-2}
\end{subfigure}\hfill%
\begin{subfigure}[t]{.3\textwidth}
	\centering
	\includegraphics[scale=\smallscale]{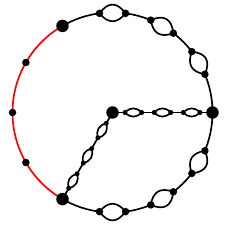}

	\smallskip
	\caption{Contracting the beads along a refined edge of $\dual(R_3(\tri))$.}
\end{subfigure}

\smallskip

\caption{The dual graphs of a 2-dimensional triangulation $\tri$ and of its 3-refinement. The four large round marks on the sub-figures b and c indicate the original vertices already present in $\dual(\tri)$.}
\label{fig:graph-refinement}
\end{figure}

To simplify the notation, we let $G = \dual(\tri)$ and $G^\ast = \dual(\tri^\ast) = \dual(R_{2\lceil \log_2 n \rceil - 1}(\tri))$. Our goal is to show that $\tww(G^\ast) \leq 4$. First we observe that it is sufficient to prove this for the case when $n = |V(G)|$ is a power of two. Indeed, if $2^{k-1} < n < 2^k$ for some positive integer $k$, we extend $\tri$ with $2^k - n$ isolated $d$-simplices, leading to a triangulation $\tri_\blacktriangle$. (Its dual graph $\dual(\tri_\blacktriangle)$ is obtained by adding $2^k - n$ isolated vertices to $G$.) Since $\log_2 |\tri_\blacktriangle| = k = \lceil \log_2 n \rceil$, the refinement parameter remains the same. Moreover, as $G^\ast$ is an induced subgraph of $\dual(R_{2\lceil \log_2 n \rceil - 1}(\tri_\blacktriangle))$ we have that $\tww(G^\ast) \leq \tww(\dual(R_{2\lceil \log_2 n \rceil - 1}(\tri_\blacktriangle))$. Henceforth, we assume that $n = 2^k$ for some positive integer $k$.

\begin{figure}[htpb]
\centering
\begin{subfigure}[t]{.48\textwidth}
	\centering
	\includegraphics[scale=\smallscale,page=2]{zipping}

	\smallskip
	\caption{The trigraph $H$, obtained by identifying the leaves of the complete, red, rooted binary tree $T_k$ of height $k$ with the set $\orig(G^\ast)$ of the original vertices of $G^\ast$. The (red) edges of $T_k$ and the refined edge $\refedge_{u,v}$ between the vertices $u,v \in \orig(G^\ast)$ are shown; the other edges of $H$ are not displayed.}
	\label{fig:zipping-refined-edge}
\end{subfigure}\hfill%
\begin{subfigure}[t]{.48\textwidth}
	\centering
	\includegraphics[scale=\smallscale,page=3]{zipping}

	\smallskip
	\caption{First, we sequentially contract the beads along the refined edge $\refedge_{u,v}$, thereby obtaining a red path $\redpath_{u,v}$ of length $2k$ between $u$ and $v$. Throughout this process, the red degree of any vertex remains at most three (realized by internal vertices of $T_k$).}
	\label{fig:zipping-red-path}
\end{subfigure}

\medskip

\begin{subfigure}[t]{.48\textwidth}
	\centering
	\includegraphics[scale=\smallscale,page=4]{zipping}

	\smallskip
	\caption{Note that there is a unique path in $T_k$ that connects $u$ and $v$. It has the same length as $\redpath_{u,v}$. In the second part, we sequentially contract pairs of vertices along these two paths that have the same distance from $u$ (and, hence, from $v$). These pairs of vertices are connected by dashed curves.}
	\label{fig:zipping-two-paths}
\end{subfigure}\hfill%
\begin{subfigure}[t]{.48\textwidth}
	\centering
	\includegraphics[scale=\smallscale,page=4]{zipping-one-contraction}

	\smallskip
	\caption{Snapshot after contracting the first pair of vertices along the paths described in the caption of \Cref{fig:zipping-two-paths}. At this point, the maximum red degree is four, and is realized by the vertex $t$ of $T_k$. As we continue ``zipping'' these two paths together, the maximum red degree remains four.}
\end{subfigure}

\medskip

\begin{subfigure}[t]{.48\textwidth}
	\centering
	\includegraphics[scale=\smallscale,page=5]{zipping}

	\smallskip
	\caption{After ``zipping together'' the two paths between $u$ and $v$, we obtain a trigraph $\widehat{H}$ that equals $H$ with the refined edge $\refedge_{u,v}$ removed. Taking $\widehat{H}$ as the starting point, we iterate the process described in \Cref{fig:zipping-refined-edge,fig:zipping-red-path,fig:zipping-two-paths} until all refined edges disappear.}
\end{subfigure}\hfill%
\begin{subfigure}[t]{.48\textwidth}
	\centering
	\includegraphics[scale=\smallscale,page=1]{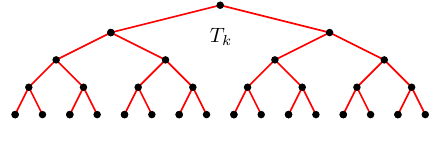}

	\smallskip
	\caption{Once all refined edges are removed, we are left with the red tree $T_k$. We conclude the process by taking an optimal contraction sequence $T_k \leadsto \bullet$, which has width three. Thus the overall maximum red degree during the entire process remains four.}
\end{subfigure}

\medskip
\caption{Illustration of a contraction sequence $H \leadsto \bullet$ of width four, showing that $\tww(H) \leq 4$.}
\label{fig:zipping}

\end{figure}

Let $\orig(G^\ast) \subset V(G^\ast)$ denote the set of original vertices of $G^\ast$. We have $|\orig(G^\ast)| = n = 2^k$. Now, let $T_k$ be the complete rooted binary tree of height $k$ with all edges already colored red. Consider the trigraph $H$ obtained by first taking the disjoint union of $G^\ast$ and $T_k$, then identifying the leaves of $T_k$ with $\orig(G^\ast)$ via an arbitrary bijection. As $G^\ast$ is a subtrigraph of $H$, we have $\tww(G^\ast) \leq \tww(H)$, hence it suffices to show that $\tww(H) \leq 4$, which can be done by exhibiting a contraction sequence $H \leadsto \bullet$ that has width 4. We use essentially the same contraction sequence as in \cite[Section 3]{berge2022deciding} with the small modification that, whenever we pick a refined edge between two original vertices $u$ and $v$ of $G^\ast$, we first contract its beads (cf.~``Additional terminology'' above) before zipping it along the unique shortest path between $u$ and $v$ in $T_k$. This process is explained and illustrated on \Cref{fig:zipping}.
\end{proof}

\subparagraph*{AI disclosure} No AI assistance was used to obtain the results presented in this paper. The second author used DeepL Write to make stylistic improvements to the text written by the authors. An autoformalized Lean proof can be found on Lax: \url{https://laxarchive.org/lax-919488}.

\newpage

\bibliography{references}

\end{document}